\documentclass[12pt]{iopart}

\usepackage{xcolor}
\usepackage{amssymb}
\usepackage{ntheorem}
\usepackage{txfonts}
\usepackage{cleveref}
\usepackage{microtype}
\usepackage{mathrsfs}
\usepackage{graphicx}
\usepackage{cite}
\usepackage{ulem}
\usepackage{caption}

\newtheorem{theorem}{Theorem}[section]
\newtheorem{proposition}{Proposition}[section]
\newtheorem{definition}{Definition}[section]
\newtheorem{lemma}{Lemma}[section]

\newtheorem{example}{Example}[section]
\newtheorem*{proof}[theorem]{Proof}
\newtheorem{remark}{Remark}[section]
\begin{document}

\title[two bifurcation sets of expansive Lorenz maps with a hole]{Two bifurcation sets of expansive Lorenz maps with a hole at the critical point}

\author{Yun Sun, Bing Li$^{\dag}$}

\address{Department of Mathematics, South China University of Technology, Guangzhou, 510461, China}
\ead{202110186724@mail.scut.edu.cn;  scbingli@scut.edu.cn}
\vspace{10pt}
\begin{indented}
\item[]April 2024
\end{indented}

\begin{abstract}
\par Let $f$ be an expansive Lorenz map on $[0,1]$ and $c$ be the critical point. The survivor set we are discussing here is denoted as $S^+_{f}(a,b):=\{x\in[0,1]:f(b)\leq f^{n}(x) \leq f(a)\ \forall n\geq0\}$, where the hole $(a,b)\subseteq [0,1]$ satisfies $a\leq c \leq b$ and $a\neq b$.
Let $a\in[0,c]$ be fixed, we mainly focus on the following two bifurcation sets:
$$ E_{f}(a):=\{b\in[c,1]:S^{+}_{f}(a,\epsilon)\neq S^{+}_{f}(a,b) \ \forall \ \epsilon>b\}, \ \ {\rm and}
$$
$$ B_{f}(a):=\{b\in[c,1]:h_{top}(S^+_{f}(a,\epsilon))\neq h_{top}(S^+_{f}(a,b)) \ \forall \ \epsilon>b\}.
$$
 By combinatorial renormalization tools, we give a complete characterization of the maximal plateau $P(b)$ such that for all $\epsilon\in P(b)$, $h_{top}(S^+_{f}(a,\epsilon))=h_{top}(S^+_{f}(a,b))$. Moreover, we obtain a sufficient and necessary condition for $E_{f}(a)=B_{f}(a)$, which partially extends the results in \cite{allaart2023} and \cite{baker2020}.

\vspace{0.5cm}
\noindent Keywords: Expansive Lorenz maps; Bifurcation sets; Survivor sets; Topological entropy; Renormalization.  \\
\noindent Mathematics Subject Classification numbers: 37E05, 37B10;

\footnotetext { Author to whom any correspondence should be addressed.}
\end{abstract}


\section{Introduction}
 \par A {\em Lorenz map} on $X=[0,1]$ is a piecewise monotone map $f:X \to X$ with a critical point $c\in (0,1)$, such that (1) $f$ is strictly increasing on $[0,c)$ and on $(c,1]$; (2) $\lim_{x \uparrow c}f(x)=1$, $\lim_{x \downarrow
	c}f(x)=0$. If, in addition, $f$ satisfies the topological expansive condition: \begin{center}
$\overline{\cup_{n \ge 0}f^{-n}(c)}=X$,
\end{center}
then $f$ is said to be an {\em expansive Lorenz map}. Lorenz maps are one-dimensional maps with a single discontinuity, which
arise as Poincar${\rm \acute{e}}$ return maps for flows on branched manifolds that model the strange
attractors of Lorenz systems. There are a lot of studies about properties of Lorenz maps, such as renormalization\cite{ding2011,glendinning1996,hubbard1990,glendinning1993}, kneading invariants\cite{cuiding2015,ds2021,glendinning1990,sun2023} and so on. For convenience, we denote by $ELM$ the set of expansive Lorenz maps, and let $LM$ be the set of Lorenz maps. Let $f\in ELM$ and $H\subset [0,1]$ be an open subinterval which is called the hole, generally, the {\em survivor set} is defined as
$$ S_{f}(H):=\{x\in X: f^{n}(x)\notin H \  \forall n\geq 0\}=X\setminus\bigcup_{n=0}^{\infty}f^{-n}(H).
$$
Notice that $S_{f}(H)$ depends on the size of $H$, the position of the hole $H$ and the map $f$.

\par There are many results concerning the Hausdorff dimension of $S_{f}(H)$. Urba${\rm\acute{n}}$ski \cite{Urbanski1986,Urbanski1987} proved that, for the case $f$ being the doubling map $T_{2}$ with a hole $(0,t)$, the dimension function $\eta_{2}:t\mapsto \dim_{\mathcal{H}}S_{2}(0,t)$ is a devil's staircase. Kalle et al.\cite{kalle2020} considered $f=T_{\beta}$ with a hole $(0,t)$, where $T_{\beta}:x\mapsto \beta x \ ( {\rm mod} \ 1)$, $x\in[0,1]$ and $\beta\in(1,2]$. They showed that the dimension function $\eta_{\beta}:t\mapsto \dim_{\mathcal{H}}S_{\beta}(0,t)$ is also a devil's staircase. Let $T_{\beta,\alpha}(x):= \beta x+\alpha \ ( \bmod \ 1)$, where $x\in[0,1]$ and $(\beta,\alpha)\in\Delta:=\{(\beta, \alpha) \in \mathbb{R}^{2}:\beta \in (1, 2], \  \alpha \in[0,2 - \beta]\}$.
 Langeveld and Samuel \cite{Langeveld2023} studied $f=T_{\beta,\alpha}$ with a hole $(0,t)$ and obtained that $\eta_{\beta,\alpha}:t\mapsto \dim_{\mathcal{H}}(S_{\beta,\alpha}(0,t))$ is a non-increasing devil's staircase. Recently, we \cite{sun2024} extended the devil staircase property to $f
 \in ELM$ with a hole $(a,b)$ at the critical point, and concerned the survivor set
  $$ S_{f}(a,b):=\{x\in [0,1]: f^{n}(x)\notin (a,b) \  \forall n\geq 0\}.$$
  We proved that if $f\in ELM$ has an ergodic absolutely continuous invariant measure, then the topological entropy function $\lambda_{f}(a):b\mapsto h_{top}(f|S_{f}(a,b))$ with $a\in[0,c]$ being fixed is a devil's staircase. Naturally, with the help of Ledrappier-Young formula due to Raith \cite{raith1994}, when considering $f=T_{\beta,\alpha}$ with a hole $(a,b)$, the Hausdorff dimension function $\eta_{f}(a):b\mapsto \dim_{\mathcal{H}}(S_{f}(a,b))$ is also a devil's staircase.

\par We emphasize that the hole $H=(a,b)$ studied in this paper satisfies $a\leq c\leq b$ and $a\neq b$. Let $f\in ELM$ with $H=(a,b)$, we always consider another definition of survivor set $S^+_f(H)$ instead of $S_f(H)$, where
$$
S^{+}_{f}(H)=\left \{
\begin{array}{ll}
\{x\in[0,1]: f^{n}(x)\geq f(b) \ \forall n\geq0\} & H=(c,b), \\
\{x\in[0,1]: f^{n}(x)\leq f(a) \ \forall n\geq0\} & H=(a,c), \\
\{x\in[0,1]:f(b)\leq f^{n}(x) \leq f(a)\ \forall n\geq0\} & H=(a,b).
\end{array}
\right.
$$
Especially, for the hole at zero,
$S^{+}_{f}(0,t)=\{x\in[0,1]: f^{n}(x)\geq t \ \forall n\geq0\}$.
 It is clear that the hole $(c,b)$ is equivalent to the hole $(0,f(b))$ for the reason that $S^{+}_{f}(c,b)=S^{+}_{f}(0,f(b))$. As a result, the hole $(a,b)$ we considered here includes the hole at zero. By \cite{sun2024}, we only need to focus on $S^{+}_{f}(H)$  since $S_{f}(H)\setminus S^+_{f}(H)$ is countable and $h_{top}(f|S_{f}(H))=h_{top}(f|S^+_{f}(H))$.
Moreover, it was also proved in \cite{sun2024} that the bifurcation set $E_f(a)$ of $S^{+}_{f}(a,b)$ is of null Lebesgue measure, where $a\in[0,c]$ is fixed and
$$ E_{f}(a):=\{b\in[c,1]:S^{+}_{f}(a,\epsilon)\neq S^{+}_{f}(a,b) \ {\rm for \ any} \ \epsilon>b \}.
$$
For any $b\in E_f(a)$, we \cite{sun2024} gave a complete characterization of the maximal plateau $I(b)$ such that for all $\epsilon\in I(b)$, $S^+_{f}(a,\epsilon)=S^+_{f}(a,b)$, and $I(b)$ may degenerate to a single point $b$. A brief introduction of the results on $S^+_{f}(H)$ is presented in Section \ref{survivor}.
\par When studying the bifurcation set of topological entropy with $a$ being fixed, we denote by
$$ B_{f}(a):=\{b\in[c,1]:h_{top}(f|S_{f}(a,\epsilon))\neq h_{top}(f|S_{f}(a,b)) \ {\rm for \ any} \ \epsilon>b\}.
$$
Example \ref{different} below shows that topological entropy may remain constant  even if the survivor set $S^{+}_{f}(a,b)$ changes. Clearly, $B_{f}(a)\subseteq E_{f}(a)$ and hence $B_f(a)$ is also of null Lebesgue measure, which is applied to prove that entropy function $\lambda_{f}(a):b\mapsto h_{top}(f|S_{f}(a,b))$ is a devil staircase in \cite{sun2024}. One question arises: can we give a full characterization of the maximal plateau $P(b)$ such that for all $\epsilon\in P(b)$, $h_{top}(f|S_{f}(a,\epsilon))=h_{top}(f|S_{f}(a,b))$? A positive answer is given in Theorem \ref{plateau}, here we start with some definitions of kneading sequences.
\par \par Let $f\in LM$ and $c$ be the critical point. The orbit of $x\in[0,1]$ by $f$ can be coded by elements of $\{0,\ 1\}^{\mathbb{N}}$. The kneading
sequence of $x$ is defined to be $\tau_f(x):=(\epsilon_0 \epsilon_1 \ldots)$, where
$$ \epsilon_i=0\ \ \ \ \ \ {\rm if}\ \ \ f^i(x)<c \ \ \ \ {\rm and} \ \ \
\epsilon_i=1\ \ \ \ \ \ {\rm if}\ \ \ f^i(x)>c.$$
This definition works for $x \notin \cup_{n \ge 0}f^{-n}(c)$. In the case where $x$ is
a preimage of $c$, $x$ has upper and lower kneading sequences $ \tau_{f}(x+)=\lim_{y\downarrow x}\tau_f(y),$ and
$\tau_{f}(x-)=\lim_{y\uparrow x}\tau_f(y),$ where the $y's$ run through
points of $[0,1]$ which are not the preimages of $c$. Given any hole $H=(a,b)$, denote by ${\rm \bf {a}}=\tau_{f}(a-)$ and ${\rm \bf {b}}=\tau_{f}(b+)$. Denote $\sigma$ as the left-shift map on $\{0,\ 1\}^{\mathbb{N}}$ and consider $$\Omega({\rm \bf {b}},{\rm \bf {a}}):=\{ \omega \in \{ 0, 1\}^{\mathbb{N}} \colon \sigma({\rm \bf {b}})\preceq  \sigma^{n}(\omega) \preceq \sigma({\rm \bf {a}}) \ \textup{for \ all} \, n \in \mathbb{N}_{0}\}.$$
By the result in {\cite[Lemma 3.4]{sun2024}}, we can always find a pair of weak-admissible (see Section \ref{lowerkneading}) kneading sequences $(1s,0t)$ such that $\Omega(1s,0t)=\Omega({\rm \bf {b}},{\rm \bf {a}})$. For simplicity, we suppose $({\rm \bf {b}},{\rm \bf {a}})$ is weak-admissible in the following theorem, and characterize the endpoints of maximal platform of identical entropy via non-periodic renormalizations (see Section \ref{determinant}). Set ${\rm \bf {b}}|_p$ as the first $p$ elements of sequence ${\rm \bf {b}}$.
\begin{theorem}\label{plateau}
Let $f\in ELM$ with $H=(a,b)$, where $a\in (0,c]$ is fixed and $b\in E_{f}(a)$. Then we have $h_{top}(f|S_{f}(a,\epsilon))=h_{top}(f|S_{f}(a,b))$ for all $\epsilon\in P(b)=[b_L,b_R]$, where
    \begin{enumerate}
\item $\tau_f(b_L+)=w_+w_-^\infty$ and $\tau_f(b_R+)=w_+^\infty$ if $({\rm \bf {b}},{\rm \bf {a}})$ can be non-periodically renormalized via $(w_+,w_-)$, including two cases $({\rm \bf {b}},{\rm \bf {a}})=(w_+^\infty,w_-w_+^\infty)$ and $(w_+w_-^\infty,w_-^\infty)$.
 \item $\tau_f(b_L+)=w_+{\rm \bf {a}}$ and $\tau_f(b_R+)=w_+^\infty$ if $({\rm \bf {b}},{\rm \bf {a}})=(w_+{\rm \bf {a}},{\rm \bf {a}})$.
\item $\tau_f(b_L+)={\rm \bf {b}}|_p{\rm \bf {a}}$ and $\tau_f(b_R+)={\rm \bf {b}}$ if ${\rm \bf {b}}$ is periodic with period $p$ and $({\rm \bf {b}},{\rm \bf {a}})$ can not be non-periodically renormalized as case (i).
\item $P(b)=\{b\}$ if ${\rm \bf {b}}$ is not periodic and $({\rm \bf {b}},{\rm \bf {a}})$ can not be non-periodically renormalized as cases (i) and (ii).
\end{enumerate}
\end{theorem}
\par We divide the proof of Theorem \ref{plateau} into two parts: the case $a\neq c$ and the case $a=c$, and prove them in a more general condition, i.e., $({\rm \bf {b}},{\rm \bf {a}})$ may not be weak-admissible. It is known that $B_f(a)\subseteq E_f(a)$, the remark below explains at which case $B_f(a)\neq E_f(a)$.
\begin{remark}\label{remark1111}
Let $a$ be fixed, $B_f(a)\neq E_f(a)$ if and only if there exists $b\in E_f(a)$ such that $I(b)\subsetneqq P(b)$.
\end{remark}

 As we can see, $I(b)\subsetneqq P(b)$ indicates that some bifurcation points of $E_f(a)$ belong to a plateau $P(b)$, which leads to $B_f(a)\subsetneqq E_f(a)$. See Example \ref{different} for an intuitive understanding. Let $a\in[0,c]$ be fixed, denote by $$ D_f(a):=\{b\in E_f(a):{\rm \bf {b}} \ \textup {is periodic}\}.$$
\begin{theorem}\label{thmxin}
Let $f\in ELM$ with a hole $(a,b)$, where $a$ is fixed. Then $D_f(a)$ is dense in $E_f(a)$.
\end{theorem}

\par  By Theorem \ref{thmxin}, for any $b\in E_f(a)$ with ${\rm \bf {b}}$ not being periodic, we can always find $b^\prime\in E_f(a)$ with ${\rm \bf {b}}^\prime$ being periodic, and the Euclidean distance between $b$ and $b^\prime$ is arbitrarily close. Notice that if $a=c$, we write the two bifurcation sets as $E_f(c)$ and $B_f(c)$. Urba${\rm\acute{n}}$ski \cite{Urbanski1986} proved that, when $f=T_2$ with a hole at zero, then $E_2(c)=B_2(c)$. Baker and Kong \cite{baker2020} considered $f=T_\beta$ ($\beta\in(1,2]$) with a hole at zero, they showed that if $\beta$ is a multinacci number, i.e., the unique root in $(1,2)$ of the equation
$$ x^{m+1}=x^m+x^{m-1}+\cdots+x+1
$$
for some $m\in \mathbb{N}$, then $E_\beta(c)=B_\beta(c)$. Recently, a characterization of $B_{\beta}(c)=E_{\beta}(c)$ for $f=T_\beta$ is given by Allaart and Kong \cite{allaart2023}. Two natural questions arise: (1) What if we study $f\in ELM$ with a hole $(a,b)$ and $a\in(0,c]$? (2) Can we give a sufficient and necessary condition for  $B_{f}(a)= E_{f}(a)$? Here we give a positive answer via linearizability condition (see definition in Subsection \ref{linearizablecondition}). Notice that for the bifurcation set $E_f(a)$ here, we only consider the case $S^+_f(H)$ being nonempty.

\begin{theorem}\label{thmplat}Let $f\in ELM$ with $H=(a,b)$, where $a\in (0,c]$ is fixed.
\begin{enumerate}
\item If $a\neq c$, then $E_{f}(a)\neq B_{f}(a)$.
\item If $a=c$, then $E_{f}(c)=B_{f}(c)$ if and only if $({\rm \bf {b}},k_-)$ is linearizable for all $b\in E_{f}(c)$.
\end{enumerate}
\end{theorem}
\begin{remark}
\
\noindent
\begin{enumerate}
\item Let $f\in ELM$ with its kneading invariants $(k_+,k_-)=(10^\infty,(01101)^\infty)$. When considering the hole $(c,b)$, we have $E_f(c)=B_f(c)$.
\item Let $f\in ELM$ with its kneading invariants $(k_+,k_-)=(10^\infty,0111000(100)^\infty)$, and the hole is $(c,b)$. Let $b$ be such that $\tau_f(b+)=(10)^\infty$, then $b\in E_f(c)$ and $((10)^\infty,k_-)$ is not linearizable. By Theorem \ref{thmplat} (ii) we have $E_f(c)\neq B_f(c)$.
\end{enumerate}
\end{remark}

\par Our work is organized as follows. In Section 2, we introduce some preliminaries, including kneading invariants, combinatorial renormalization and kneading determinants. Section 3 gives some valuable lemmas and presents the proof of Theorem \ref{plateau}. In Section 4, we obtain a sufficient and necessary condition for $E_f(a)=B_f(a)$.

\section{Preliminaries}

\subsection{Kneading invariants}\label{lowerkneading}
We equip the space $\{0,1\}^\mathbb{N}$ of infinite sequences with the topology induced by the usual metric $d \colon \{0,1\}^\mathbb{N} \times \{0,1\}^\mathbb{N} \to \mathbb{R}$ which is given by
$$
\label{cases}
d(\omega, \nu):=\cases{0& if $\omega = \nu$,\\
2^{- |\omega \wedge \nu| + 1}& otherwise.\\}
$$
Here $|\omega \wedge \nu | := \min \, \{ \, n \in \mathbb{N} \colon \omega_{n} \neq \nu_n \}$, for all $\omega = (\omega_{1}\omega_{2}\dots) , \nu = ( \nu_{1} \nu_{2}\dots) \in \{0, 1\}^{\mathbb{N}}$ with $\omega \neq \nu$, and we denote the lexicographic order as $\omega\prec \nu$ if $\omega_n<\nu_n$.
Note that the topology induced by $d$ on $\{ 0, 1\}^{\mathbb{N}}$ coincides with the product topology on $\{ 0, 1\}^{\mathbb{N}}$.  For $n\in\mathbb{N}$ and $\omega\in\{0,1\}^{\mathbb{N}}$, we set $\omega|_{1}^{n}=\omega|_{n}=(\omega_{1} \dots \omega_{n})$ and call $n$ the length of $\omega|_{n}$. We denote $\sigma \colon \{ 0, 1 \}^{\mathbb{N}} \circlearrowleft$ as the \textsl{left-shift map} which is defined by $\sigma(\omega_{1} \omega_{2} \dots) \coloneqq (\omega_{2} \omega_{3}\dots)$.  A \textsl{subshift} is any closed subset $\Omega \subseteq \{0,1\}^\mathbb{N}$ such that $\sigma(\Omega) \subseteq \Omega$.  Given a subshift $\Omega$ and $n \in \mathbb{N}$ we set $
\Omega|_{n} := \left\{ (\omega_{1} \dots \omega_{n}) \in \{ 0, 1\}^{n} \colon \textup{there \ exists} \ \omega \in \Omega \, \ \textup{with} \ \, \omega|_{n} = (\omega_{1} \dots \omega_{n})\right\}
$
and denote by $\Omega^{*}:= \bigcup_{n = 1}^{\infty} \Omega|_{n}$ for the collection of all finite words. For $\xi\in\Omega^{*}$, we denote $|\xi|$ as the length of $\xi$, and $\#\Omega|_{n}$ as the cardinality of $\Omega|_{n}$. Moreover, we write $\omega = (\omega_{1} \omega_{2} \dots \omega_{n})^\infty $ if $\omega$ is periodic with period $n$, and $\omega = \omega_{1}\dots \omega_{k} (\omega_{k+1} \dots \omega_{k+n})^{\infty}$ if $\omega$ is eventually periodic.

\par We call an infinite sequence $\omega\in\{0,1\}^{\mathbb{N}}$ is {\bf self-admissible} if $\sigma(\omega)$ is lexicographically largest or lexicographically smallest, that is, $\sigma (\omega)\succeq\sigma^{n}(\omega)$ for all $n\geq0$, or $\sigma (\omega)\preceq\sigma^{n}(\omega)$ for all $n\geq0$. Let $f\in ELM$ with a hole $(a,b)$, and denote by ${\rm \bf {a}}=\tau_{f}(a-)$ and ${\rm \bf {b}}=\tau_{f}(b+)$. By the proof of {\cite[Lemma 3.1]{sun2024}}, it is trivial that $S_{f}(a,b)\subseteq\{0,1\}$ if ${\rm \bf {a}}|_2=00$ or ${\rm \bf {b}}|_2=11$. Hence we only consider the case that ${\rm \bf {a}}|_2=01$ and ${\rm \bf {b}}|_2=10$ here.
 Observed that sometimes ${\rm \bf {a}}$ or ${\rm \bf {b}}$ may not be self-admissible, for instance, there exists an integer $k$ such that $\sigma ({\rm \bf {a}})\prec\sigma^{k}({\rm \bf {a}})$, but we can verify that $({\rm \bf {a}}|_k)^\infty$ is self-admissible. As a result, we can always obtain two self-admissible sequences for a given hole $(a,b)$, hence ${\rm \bf {a}}$ and ${\rm \bf {b}}$ are considered to be self-admissible throughout this paper.

\par Let $f\in LM$,
$(k_+, k_-)=(\tau_f(c+), \tau_f(c-)) $ are called the {\bf kneading
invariants} of $f$, which were used to developing the combinatorial
theory of expansive Lorenz map.
The kneading space of $f$, also called Lorenz-shift, is $\Omega(f)=\{\tau_f(x): x \in I\}.$
Since $\sigma$ is the shift map operating on the Lorenz shift
$\Omega(f)$, then clearly $\tau_{f}(f(x))=\sigma(\tau_f(x)))$, with similar
results holding for the upper and lower kneading sequences of points
$x$ which are pre-images of $c$. The dynamics of $f$ on
$I$ can be modeled by $(\sigma,\Omega(f))$.

\begin{theorem}{\bf\textup{({\cite[Theorem 2]{hubbard1990}})}} \label{rem:k+k-}
Let $f$ be a Lorenz map. Then the kneading space $\Omega(f)$ is completely determined by the kneading invariants $(k_+,k_-)$ of $f$;  indeed, we have
\begin{eqnarray*}
\Omega(f) = \Omega(k_+,k_-) =\left\{ \omega \in \{ 0, 1\}^{\mathbb{N}} \colon \sigma(k_{+})\preceq  \sigma^{n}(\omega) \preceq \sigma(k_{-}) \ \textup{for \ all} \, n \in \mathbb{N}_{0} \right\}.
\end{eqnarray*}
\end{theorem}

\par Clearly, $\Omega(f)$ is closed with respect to the metric $d$ and hence is a subshift, and both $k_+$ and $k_-$ are self-admissible. We see that each $f\in ELM$ corresponds to a pair of kneading invariants, however, what kind of sequences in $\{ 0, 1 \}^{\mathbb{N}}$ can be the kneading invariants of $f\in ELM$? It was stated by Hubbard and Sparrow \cite{hubbard1990} as follows and we call it {\bf H-S admissibility condition}.

\begin{theorem}{\bf\textup{({\cite[Theorem 1]{hubbard1990}})}} \label{kneading space}
If $f\in ELM$, then its kneading invariants $(k_{+},k_{-})$ satisfy
	\begin{equation} \label{expanding}
		\sigma (k_{+})\preceq \sigma^n(k_{+})\prec\sigma (k_{-}), \ \ \ \ \ \ \sigma (k_{+})\prec \sigma^n(k_{-})\preceq\sigma (k_{-}) \ \ \ \  \forall n \ge 0,
	\end{equation}
 Conversely, given any two sequences $k_{+}$ and $k_{-}$ satisfying $(\ref{expanding})$, there exists an $f\in ELM$ with $(k_{+},k_{-})$ as its kneading invariant, and $f$ is unique up to conjugacy.
\end{theorem}

\par In addition, we give the definition of {\bf weak-admissible}, which also includes the non-expansive cases. We say the kneading invariants $(k_{+},k_{-})$ are weak-admissible if satisfying
	\begin{equation} \label{non-expanding}
\sigma(k_{+}) \preceq\sigma^{n}(k_{+})\preceq\sigma(k_{-}) \ {\rm and} \ \sigma(k_{+})\preceq\sigma^{n}(k_{-})\preceq\sigma(k_{-}) \ \ \ \ {\rm for  \ all}\ \ n\geq 0,
	\end{equation}
which means there may exist $n$ such that $\sigma^{n}(k_{+})=k_{-}$ or $\sigma^{n}(k_{-})=k_{+}$. Clearly, H-S admissibility implies weak admissibility, but not vice versa. There are many trivial cases of weak admissible kneading invariants, such as the kneading invariants induced by rational rotations, and some renormalizable (see the next subsection) Lorenz maps.

\subsection{Combinatorial renormalization}\label{determinant}
Renormalization is a central concept in contemporary dynamics. The idea of renormalization for Lorenz map was
introduced in studying simplified models of Lorenz attractor,
apparently firstly in Palmer \cite{palmer1979} and Parry \cite{parry1979lorenz}.  Here we focus on the renormalization in combinatorial way. Glendinning
and Sparrow \cite{glendinning1993} presented a comprehensive study of the
renormalization by investigating the kneading invariants of
expanding Lorenz maps. The following definition is essentially from \cite{glendinning1993}.


\begin{definition} \rm
Let $f$ be a Lorenz map, we say the kneading invariants $K=(k_{+},k_{-})$ of $\Omega(f)$ is \textsl{renormalizable} if there exist finite, non-empty words $(w_+,\ w_-)$, such that
\begin{equation} \label{com-renormal1}
\left \{ \begin{array}{ll}
k_{+} =&w_+ w_-^{p_1} w_+^{p_2} \cdots,\\
k_{-} =&w_- w_+^{m_1} w_-^{m_2} \cdots,
\end{array}
\right.
\end{equation}
where $w_+=10\cdots$, $w_-=01\cdots$, the lengths $|w_+|>1$
and $|w_-|> 1$, and $p_1, m_1>0$. The kneading invariants of the \textsl{renormalization} is $RK=(R^{1}k_+, R^{1}k_-)$, where
\begin{equation} \label{com-renormal2}
\left \{ \begin{array}{ll}
R^{1}k_+=&1 0^{p_1} 1^{p_2} \cdots,\\
R^{1}k_-= & 0 1^{m_1} 0^{m_2} \cdots.
\end{array}
\right.
\end{equation}
\end{definition}

\par To describe the renormalization more concisely, we use the $*$-product, which is introduced in \cite{glendinning1990}. The $*$-product of renormalizaition is defined to be $K=(k_+, k_-)=W*R^1K$, i.e.,
\begin{equation}\label{product}
(k_{+},k_{-})=(w_+,\ w_-)*(R^{1}k_+,R^{1}k_-)
\end{equation}
where $(k_{+},k_{-})$ is the pair of sequences obtained by
replacing $1$'s by $w_+$, replacing $0$'s by $w_-$ in $R^{1}k_+$ and $R^{1}k_-$.
Using $*$-product, (\ref{com-renormal1}) and (\ref{com-renormal2}) can be expressed by (\ref{product}). So $(k_+,k_-)$ is renormalizable if and only if it can be decomposed as the $*$-product; otherwise, we say that $(k_{+},k_{-})$ is {\bf prime}. One can check that the $*$-product satisfies the associative law, but does not satisfy the commutative law in general. Note that we do not involve $(w_+, w_-)=(1, 01)$ and $(w_+, w_-)=(10, 0)$ in the definition of renormalization as in \cite{glendinning1993}, since these two cases correspond to trivial renormalizations. Furthermore, $R^1K=W^{\prime}\ast R^{2}K=(w^{\prime}_+,\ w^{\prime}_-)*(R^{2}k_+,R^{2}k_-)$ if $R^1K$ is also renormalizable. And $(k_+,k_-)$ is $m$ $(0\leq m\leq\infty)$ times renormalizable if the renormalization process can proceed $m$ times. We denote by $(R^{m}k_+,R^{m}k_-)$ the kneading invariants after $m$ times renormalizations.

\par Before introducing periodic renormalization, we first see some properties about rational rotations, that is, $T_{\beta,\alpha}$ with $\beta=1$ and $\alpha\in\mathbb{Q}\cap(0,1)$. Similarly, we can also obtain the upper and lower kneading sequences of the critical point $c_{1,\alpha}=1-\alpha$, which correspond to the kneading invariants $(k_{+},k_{-})$ of $T_{1,\alpha}$. For any $\alpha=p/q\in  \mathbb{Q}\cap(0,1)$, we know that the orbits $\{T_{1,\alpha}(x)\}$ are periodic with period $q$ for all $x\in[0,1]$, and it is easy to obtain that $T_{1,\alpha}(0)=T_{1,\alpha}(1)$. As a result, if we denote $k_{+}=v^{\infty}=(v_{1}\cdots v_{q})^{\infty}$ and $k_{-}=u^{\infty}=(u_{1}\cdots u_{q})^{\infty}$, then we have the following two properties about the finite words $u$ and $v$,
\par (1) $k_{+}$ and $k_{-}$ are on the same periodic orbit, i.e., there exists an integer $s\leq q-1$ such that $\sigma^{s}(k_{+})=k_{-}$.
\par (2)  $v|_{2}=10$, $u|_{2}=01$, and $v_{i}=u_{i}$ for all $i\in\{3,\cdots,q\}$.\\
For convenience, we call such two finite words $(v,u)$ be {\bf rational words} since they correspond to a rational rotation.

\par
A renormalization is called {\bf periodic renormalization} if the finite renormalization words $w_{+}$ and $w_{-}$ are rational words; otherwise, it is called {\bf non-periodic renormalization}. As we can see, if kneading invariants $(k_{+},k_{-})$ can be periodically renormalized via renormalization words $w_{+}$ and $w_{-}$, then there exists a rational number $\alpha=p/q\in\mathbb{Q}\cap(0,1)$ such that $(w_{+}^{\infty}, w_{-}^{\infty})$ are the kneading invariants of rational rotation $T_{1,\alpha}$. Hence it is well understood that periodic renormalization is also called primary $q(p)$-cycle in \cite{glendinning1990,palmer1979}. See the following example.
\begin{example}\label{periodic}(Periodic renormalization)\rm
\begin{enumerate}
\item Let $(k_{+},k_-)=((100101)^{\infty},(0110)^{\infty})$, then $(k_{+},k_{-})$ can be renormalized to $RK=((100)^{\infty},(01)^{\infty})$ via words $w_{+}=10$ and $w_{-}=01$. Here is periodic renormalization since $(w_{+}^{\infty}, w_{-}^{\infty})$ corresponds to the kneading invariants of rational rotation $T_{1,1/2}$. It is clear that $h_{top}(\sigma,\Omega(k_+,k_-))>0$ and $(k_{+},k_{-})$ can only be periodically renormalized for finitely many times, hence $(k_{+},k_{-})$ is also linearizable.
\item Let $(k_{+},k_-)=((100011)^{\infty},(011100)^{\infty})$, then $(k_{+},k_{-})$ can be renormalized via words $w_{+}=100$ and $w_{-}=011$. However, $(w_{+}^{\infty}, w_{-}^{\infty})$ does not correspond to the kneading invariants of any rational rotation, hence it is non-periodic renormalization. Moreover,  $(k_{+},k_{-})$ does not satisfy Definition \ref{defn:admissible} (2), hence $(k_{+},k_{-})$ is non-linearizable.
\end{enumerate}
\end{example}

\par Let $f\in ELM$ with a hole $(a,b)$. Recall that
$$\Omega({\rm \bf {b}},{\rm \bf {a}})=\{ \omega \in \{ 0, 1\}^{\mathbb{N}} \colon \sigma({\rm \bf {b}})\preceq  \sigma^{n}(\omega) \preceq \sigma({\rm \bf {a}}) \ \textup{for \ all} \, n \in \mathbb{N}_{0}\}.$$
We say that $\Omega({\rm \bf {b}},{\rm \bf {a}})$ is renormalizable if the pair of kneading sequences $(1s,0t)$ is renormalizable, where $\Omega({\rm \bf {b}},{\rm \bf {a}})=\Omega(1s,0t)$ and $(1s,0t)$ is weak-admissible. For instance, let $(k_+,k_-)=((10011)^\infty,(0111010100)^\infty)$, by Lemma \ref{exist weak-admissible}, we have $1s=(10011)^\infty$ and $0t=(0111010)^\infty$. Hence we say $\Omega(k_+,k_-)$ can be non-periodically renormalized via $w_+=10$ and $w_-=011$.

\subsection{Linearizable kneading invariants}\label{linearizablecondition}

As one of the simplest piecewise linear maps on the interval, the linear mod one transformation defined by $T_{\beta, \alpha}(x)=\beta x+\alpha \  \ ( \bmod \ 1)$
has attracted considerable attention.
It was proved that $T_{\beta,\alpha}$ has a unique absolutely continuous invariant measure which is ergodic, and this measure is the unique measure of maximal entropy with entropy $\log \beta$. The standard definition for the topological entropy of continuous maps using $(n,\epsilon)$-separated sets can be used to define entropy for piecewise continuous maps. An alternative way of calculating entropy, which is particularly convenient for the symbolic approach, is via the cardinality of finite words, i.e.,
$$ h_{top}(f)=h_{top}(\sigma,\Omega(f))=\displaystyle{\lim_{n \to \infty}  \frac{\log \left( \# \Omega(f)|_{n}  \right)}{n}}.
$$
The limit exists for the sequence $\{ \log(\# \Omega(f)\vert_{n} \}_{n \in \mathbb{N}}$ is sub-additive, and hence
$$ h_{top}(f)=\displaystyle{\inf_{n}  \frac{\log \left( \# \Omega(f)|_{n}  \right)}{n}}.
$$

\par Since $T_{\beta, \alpha}$ is a uniformly linear Lorenz map, a natural question arises that, when is an expansive Lorenz map uniformly linearizable, that is, topologically conjugate to some $T_{\beta, \alpha}$?  Denote $\Omega_{\beta,\alpha}$ be the $\beta$-shift induced by $T_{\beta,\alpha}$.

\begin{lemma}[{\cite[Theorem A]{glendinning1990}}{\cite[Theorem A]{cuiding2015}}]\label{finite-periodic}
When $\sqrt{2}<\beta<2$, for all $\alpha\in(0,2-\beta)$, $\Omega_{\beta,\alpha}$ is prime. When  $1<\beta\leq\sqrt{2}$ and $\alpha\in(0,2-\beta)$, $\Omega_{\beta,\alpha}$ is either prime or can only be periodically renormalized for finitely many times.
\end{lemma}

\par By Lemma \ref{finite-periodic} and Theorem \ref{kneading space}, we are able to give the definition of linearizable kneading invariants for $f\in ELM$.
Since here we always consider the $\beta$-shifts with $\beta>1$, an extra condition is needed to make sure the topological entropy is bigger than $0$. Recall that $\Omega(k_{+}, k_{-}):=\{\omega\in\{0,1\}^{\mathbb{N}}: \sigma(k_{+})\preceq\sigma^{n}(\omega) \preceq \sigma(k_{-}) \ \textup{for all} \, n \in \mathbb{N}_{0}\}$.

\begin{definition}\label{defn:admissible} \rm
A pair of infinite sequences $(k_{+},k_{-})$ is said to be \textsl{linearizable} if the following conditions are satisfied:
\begin{enumerate}
\item\label{enumi2:defn_admissible} 
$(k_{+},k_{-})$ is H-S admissible,
\item\label{enumi4:defn_admissible} $(k_{+},k_{-})$ is prime or can only be periodically renormalized for finitely many times,

\item\label{enumi3:defn_admissible} $h_{top}(\sigma,\Omega(k_{+}, k_{-})) > 0$.
\end{enumerate}
\end{definition}

\par  In conclusion, two infinite sequences $k_{+},k_{-}\in \{ 0, 1\}^{\mathbb{N}}$ are kneading invariants for an intermediate $\beta$-shift if and only if $(k_{+},k_{-})$ is linearizable, i.e., satisfying Definition \ref{defn:admissible}. For the intuitive explanation of linearizable and non-linearizable cases, also see Example \ref{periodic}.

\subsection{Kneading determinant and entropy}\label{determinant2}

The ideas for kneading determinant goes back to \cite{milnor1988}; see also \cite{glendinning1996}. Let $({k_ + },{k_ - })$ be the kneading invariants of $f\in LM$, where $k_{+}=(v_{1}v_{2}\cdots)$ and $k_{-}=(w_{1}w_{2}\cdots)$. Then the kneading determinant is a formal power series defined as $K(t) = {K_ + }(t) - {K_ - }(t)$ , where
$${K_ + }(t) = \sum\limits_{i = 1}^\infty  {v_{i}{t^{i-1}}},\ \ \ \ {K_ - }(t) = \sum\limits_{i = 1}^\infty  {w_{i}{t^{i-1}}}.$$
 The following lemma offers a straight method to calculate $h_{top}(f)$ if its kneading invariants $(k_+,k_-)$ are given.

\begin{lemma}[{\cite[Theorem 3]{glendinning1996}}{\cite[Lemma 3]{barnsley2014}}]\label{zeros}
Let $({k_ + },{k_ - })$ be the kneading invariants of $f\in LM$ with $h_{top}(f)>0$, and $K(t)$ be the corresponding kneading determinant. Denote $t_{0}$ be the smallest positive root of $K(t)$ in $(0,1)$, then $h_{top}(f)=-\log t_{0}$.
\end{lemma}

Naturally, if $(k_{+},k_{-})$ is a pair of linearizable kneading invariants, i.e, $(k_{+},k_{-})$ corresponds to an intermediate $\beta$-shift, then $1/\beta$ equals the smallest positive root of $K(t)$ in $(0,1)$. When $(k_+,k_-)$ can be non-periodically renormalized, some interesting phenomena related with entropy will appear, see the following lemma.

\begin{lemma}[{\cite[Lemma 8]{barnsley2014}}{\cite[Proposition 2]{ds2021}}]\label{xunibeta}
Let $({k_ + },{k_ - })$ be the kneading invariants of $f\in LM$. If $(k_{+},k_{-})$ can be non-periodically renormalized via renormalization words $(w_{+}, w_{-})$, then $h_{top}(f)=h_{top}(\sigma,\Omega(w_{+}^\infty,w_{-}^\infty))$.
\end{lemma}
\par We call $f\in LM$ is non-expansive if $f$ has homtervals $J$ on which the kneading sequence is constant (i.e., $f^{n}|_{J}$ is a homeomorphism for all $n\geq 0$). What follows are some results from the proof of \textup{{\cite[Theorem 3]{glendinning1996}}}, which is used to calculate the entropy for the non-expansive cases.

\begin{lemma}{\bf\textup{({\cite[Theorem 3]{hubbard1990}})}} \label{normalrenor}
Let $f\in LM$ with its kneading invariants $(k_+,k_-)$. If $(k_+,k_-)$ can be renormalized via $(w_+,w_-)$, then
$$ h_{top}(\Omega(w_+w_-^{\infty},w_-^{\infty}))=h_{top}(\Omega(w_+^{\infty},w_-w_+^{\infty}))
=h_{top}(\Omega(w_+^{\infty},w_-^{\infty})),
$$
and
$$ h_{top}(\Omega(w_+k_-,k_-))=
h_{top}(\Omega(w_+^{\infty},k_-)), \ h_{top}(\Omega(k_+,w_-k_+))=
h_{top}(\Omega(k_+,w_-^{\infty})).
$$

\end{lemma}

\par Using Lemma \ref{xunibeta} and Lemma \ref{normalrenor}, we are able to construct lots of Lorenz maps with the same entropy via non-periodic renormalization, and we can even let them have the same upper kneading invariant $k_{+}$ or the same lower kneading invariant $k_-$. See the following Example \ref{different} for an intuitive understanding of this construction. As an application in open dynamical systems, we obtain the maximal plateau of the same entropy, which is stated as Theorem \ref{plateau}. For the entropy related with periodic renormalizations, we have the following lemma from \cite{ds2021}.

\begin{example}\label{different} (Kneading invariants with same entropy)\rm
\
\par Here we construct three differen Lorenz maps with the same lower kneading sequence $k_{-}$ and meanwhile with the same entropy by non-periodic renormalization, actually there are uncountably many such Lorenz maps. Let $k_{-}=(01110)^{\infty}$ be fixed, and consider the following three different Lorenz maps with different upper kneading sequences,
\begin{equation*}
\left \{ \begin{array}{ll}
f:  \ (k_{+}^{1},k_-)=((10011)^{\infty},(01110)^{\infty}),\\
g:  \ (k_{+}^{2},k_-)=((10011011)^{\infty},(01110)^{\infty}),\\
h:  \ (k_{+}^{3},k_-)=(10(011)^{\infty},(01110)^{\infty}).
\end{array}
\right.
\end{equation*}
We can see that all of them can be non-periodically renormalized via $w_{+}=10$ and $w_{-}=011$, by Lemma \ref{xunibeta} above, they all have the same entropy with $h_{top}(\Omega((10)^{\infty},(011)^{\infty}) )=1.3247$. As a result, although their kneading spaces $\Omega(f)\neq\Omega(g)\neq\Omega(h)$, they have the same entropy and the same lower kneading invariant.
\end{example}

\begin{lemma} {\bf\textup{({\cite[Proposition 3.5]{ds2021}})}}\label{periodicentropy}
Let $f$ be a linearizable expansive Lorenz map with $m$ $(0\leq m<\infty)$ times periodic renormalizations and $(k_+,k_-)$ be its kneading invariants. Then $$h_{top}(f)=\frac{\log\beta}{l_{1}l_{2}\ldots l_{m}},$$
where $1/\beta$ is the smallest positive root of kneading determinant induced by $(R^{m}k_+,R^{m}k_-)$, and $l_{i}$ means the length of the $i$th periodic renormalization words.
\end{lemma}

\begin{remark}\label{atleastoncerenor}
Let $f\in ELM$ be $m$ $(0\leq m\leq\infty)$ times renormalizable and $j$-th $(1\leq j<\infty)$ renormalization be the nearest non-periodic renormalization with words $(w_+,w_-)$. Then all the $i$-th ($0\leq i\leq j-1$) renormalizations are periodic, and denote their renormalization words as $(w_{+i},w_{-i})$ and $l_i =|w_{+i}|=|w_{-i}|$. Using $\ast$-product, we have
 $$(k_+,k_-)=(w_{+1},w_{-1})\ast \cdots \ast(w_{+(j-1)},w_{-(j-1)})\ast(R^{j-1}k_+,R^{j-1}k_-).$$
Denote $h_{top}(\Omega(w_{+}^\infty,w_{-}^\infty))=\log \beta$. Applying Lemma \ref{xunibeta} and Lemma \ref{periodicentropy}, we have
$$ h_{top}(\Omega(R^{j-1}k_+,R^{j-1}k_-))=h_{top}(\Omega(w_{+}^\infty,w_{-}^\infty))=\log \beta.
$$
and
$$h_{top}(\Omega(k_+,k_-))=\frac{\log\beta}{l_1 l_2 \cdots l_{j-1}}.$$
\end{remark}

\subsection{Results on survivor set $S^+_{f}(a,b)$}\label{survivor}

\par Denote $\tilde{S}^+_{f}(H)$ as the symbolic representation of $S^+_{f}(H)$. To facilitate the proof of our results, we list the following results about the survivor set $S^{+}_{f}(a,b)$, which are essentially from \cite{sun2024}.

\begin{proposition}  {\bf\textup{({\cite[Theorem 1.1, Proposition 4.3]{sun2024}})}}\label{minuscountable}
The set $S_{f}(a,b)\setminus S^+_{f}(a,b)$ is countable and $h_{top}(f|S_{f}(a,b))=h_{top}(f|S^+_{f}(a,b))$.
\end{proposition}

\begin{lemma} {\bf\textup{({\cite[Lemma 3.4]{sun2024}})}}\label{exist weak-admissible}
If $\tilde{S}^{+}_{f}(a,b)\nsubseteqq\{0^\infty,1^\infty\}$ and $({\rm \bf {b}},{\rm \bf {a}})$ is not weak-admissible, then there exist weak-admissible kneading sequences $(1s,0t)$ such that $\Omega(1s,0t)=\tilde{S}^+_{f}(a,b)$.
\end{lemma}

\begin{theorem} {\bf\textup{({\cite[Theorem 1.4]{sun2024}})}}\label{plateauxin}
Let $f\in ELM$ with a hole $(a,b)$. If $\rm \bf {b}$ is periodic, then there exists a maximal plateau $I(b)$ such that for all $\epsilon\in I(b)$, $S^{+}_{f}(a,\epsilon)=S^{+}_{f}(a,b)$. The endpoints of $I(b)$ are also characterized.
\end{theorem}

\begin{proposition} {\bf\textup{({\cite[Proposition 4.1]{sun2024}})}}\label{bifurcationsetEa}
We have $E_{f}(a)=\{b\in[c,1]: {\rm \bf {b}}\in \tilde{S}^{+}_{f}(a,b)\}=\{b\in[c,1]:\sigma({\rm \bf {b}}) \preceq  \sigma^{n}({\rm \bf {b}}) \preceq \sigma({\rm \bf {a}}) \ \forall \ n\geq0\}$, where $a$ is fixed, and $E_{f}(a)$ is of null Lebesgue measure.
\end{proposition}
\begin{theorem}{\bf\textup{({\cite[Theorem 1.6]{sun2024}})}}\label{devil}
Let $f\in ELM$ with ergodic a.c.i.m.. Then the topological entropy function $\lambda_{f}(a):b\mapsto h_{top}(\tilde{S}_{f}(a,b))$ is a devil staircase, where $a$ is fixed.
\end{theorem}
By the results from \cite{sun2024}, we conclude the following remark on bifurcation points in $E_{f}(a)$.
\begin{remark}\label{aperiodic11}
Let $b\in E_{f}(a)$. If ${\rm \bf {b}}$ is periodic, then $I(b)$ is a subinterval of $[c,1]$ and $b$ is the right endpoint. If ${\rm \bf {b}}$ is not periodic, then $I(b)=\{b\}$.
\end{remark}

\section{Plateau of entropy}\label{platformsec}
A complete characterization of the plateau of the survivor set $\tilde{S}^{+}_{f}(H)$ was given by Theorem \ref{plateauxin} above. It can be seen from Example \ref{different} that even when the survivor set $\tilde{S}^{+}_{f}(H)$ changes, the topological entropy may remain unchanged. A natural question arises: Can we give a full characterization of the plateau of entropy $h_{top}(S_{f}(H))$? We obtain a positive answer in Theorem \ref{plateau} and prove it in different cases. Here we only consider the plateaus with positive entropy.

\begin{lemma} \label{w+}
Let $(w_+,w_-)$ be non-periodic renormalization words and the lower kneading sequence $k_-$ be fixed. If $k_-$ consists of $(w_+,w_-)$, then
$$ h_{top}(\Omega(w_+^{\infty},k_-))=h_{top}(\Omega(w_+w_-^{\infty},k_-))
=h_{top}(\Omega(w_+^{\infty},w_-^{\infty}))>0.
$$
\begin{proof}\rm
 We know that the lower kneading sequence $k_-$ of a Lorenz map must be self-admissible. If $k_-$ consists of $(w_+,w_-)$, it is clear that $(w_+w_-^{\infty},k_-)$ is weak-admissible kneading invariants. By Lemma \ref{xunibeta}, $\Omega(w_+w_-^{\infty},k_-)$ can be non-periodically renormalized via $(w_+,w_-)$ and naturally
$h_{top}(\Omega(w_+w_-^{\infty},k_-))
=h_{top}(\Omega(w_+^{\infty},w_-^{\infty}))$. On the other hand, $(w_+^{\infty},k_-)$ is weak-admissible if and only if $k_-=w_-w_+^\infty$, and $h_{top}(\Omega(w_+^{\infty},w_-w_+^\infty))=h_{top}(\Omega(w_+^{\infty},w_-^{\infty}))$ by Lemma \ref{normalrenor}. At the case $k_-\neq w_-w_+^\infty$, it can be seen that $(w_+^{\infty},k_-)$ is not weak-admissible. Using Lemma \ref{exist weak-admissible}, we can find a pair of weak-admissible kneading invariants $(w_+^{\infty},w_-^{\infty})$ such that $\Omega(w_+^{\infty},w_-^{\infty})=\Omega(w_+^{\infty},k_-)$, hence they have the same entropy. Moreover, by the assumption that $k_-$ consists of $(w_+,w_-)$, $k_+=w_+^\infty$ and $k_+=w_+w_-^{\infty}$ are the lexicographically largest and smallest sequence consisting of $(w_+,w_-)$, respectively, such that $\Omega(k_+,k_-)$ has the same entropy as $h_{top}(\Omega(w_+^{\infty},w_-^{\infty}))$. $\hfill\square$
\end{proof}
\end{lemma}
\par For the simplicity of our proof, here we give a wider definition of non-periodic renormalization based on Lemma \ref{w+}. From now on, both $(w_+w_-^{\infty},w_-^\infty)$ and $(w_+^{\infty},w_-w_+^\infty)$ can also be regarded as non-periodic renormalizable kneading invariants.


\begin{lemma}\label{lemmaxin222}
Let $(k_+,k_-)$ be H-S admissible and $(1s,0t)$ be  weak-admissible with $\Omega(k_+,k_-)\neq\Omega(1s,0t)$, and they have positive entropy. If neither $k_+$ nor $k_-$ is periodic and $(k_+,k_-)$ can not be non-periodically renormalized, then we have $h_{top}(\Omega(k_+,k_-))\neq h_{top}(\Omega(1s,0t))$.
\begin{proof} \rm
Without loss of generality, here we only consider $\Omega(k_+,k_-)\subsetneqq\Omega(1s,0t)$, i.e., $ k_-\preceq0t\prec1s\preceq k_+ $. By the monotonicity of topological entropy, we have $h_{top}(\Omega(k_+,k_-))\leq h_{top}(\Omega(1s,k_-))$, next we prove $h_{top}(\Omega(k_+,k_-))\neq h_{top}(\Omega(1s,k_-))$ by contradiction.
H-S admissibility of $(k_+,k_-)$ indicates that two weak-admissible cases $(w_+k_-,k_-)$ and $(k_+,w_-k_+)$ are excluded. Hence $(k_+,k_-)$ is linearizable, then there exists $(\beta,\alpha)\in \Delta$ such that $\Omega_{\beta,\alpha}=\Omega(k_+,k_-)$, where $\beta\in[1,2)$ and $\alpha\in[0,2-\beta]$. Although $(1s,0t)$ may not be linearizable, we can still calculate a pair of $(\beta^\prime,\alpha^\prime)$  by the following formula in \cite{ds2021},
\begin{equation}\label{gongshi1} \alpha  = (\beta  - 1)(\sum\limits_{i = 0}^\infty  {\frac{{{k_ + }(i)}}{{{\beta ^i}}} - 1)} = (\beta  - 1)(\sum\limits_{i = 0}^\infty  {\frac{{{k_ - }(i)}}{{{\beta ^i}}} - 1),}
\end{equation}
where $k_ \pm(i)$ means the $i$-th symbol of $k_\pm$. They have the same entropy implies that $\beta=\beta^\prime$, and $\alpha=\alpha^\prime$ by $\alpha\leq\alpha^\prime$ and $\alpha\geq\alpha^\prime$.
\par Case 1, $(1s,0t)$ is also linearizable. Then $(\beta,\alpha)=(\beta^\prime,\alpha^\prime)$ indicates that $(k_+,k_-)=(1s,0t)$, which contradicts with $\Omega(k_+,k_-)\subsetneqq \Omega(1s,0t)$.
\par Case 2, $(1s,0t)$ is not linearizable. If $(1s,0t)$ can be non-periodically renormalized via $(w_+,w_-)$, including the cases $(1s,0t)=(w_+,w_-w_+^\infty)$ and $(w_+w_-^\infty,w_-^\infty)$. By {\cite[Lemma 8]{barnsley2014}} and {\cite[Proposition 2]{ds2021}},
 $(\beta,\alpha)=(\beta^\prime,\alpha^\prime)$ indicates that $(k_+,k_-)=(w_+^\infty,w_-^\infty)$, which contradicts with our condition. If  $(1s,0t)$ is weak-admissible, that is, $(1s,0t)=(1s,w_-1s)$ (or $(w_+0t,0t)$), we have $h_{top}(\Omega(1s,0t))=h_{top}(\Omega(1s,w_-^\infty))$ ( or $h_{top}(\Omega(w_+^\infty,0t)$ ) by Lemma \ref{normalrenor}. Similarly, we have that $(k_+,k_-)=(1s,w_-^\infty)$ (or $(w_+^\infty,0t)$), which contradicts with our assumption that neither $k_+$ nor $k_-$ is periodic.
$\hfill\square$
\end{proof}
\end{lemma}
\begin{remark}\label{remarkxin}

Let $(k_+,k_-)$ be H-S admissible and $(1s,0t)$ be  weak-admissible with $\Omega(k_+,k_-)\neq\Omega(1s,0t)$ and $1s\succ k_+$, and they have positive entropy. If $k_+$ is not periodic, $(k_+,k_-)$ can not be non-periodically renormalized and $(k_+,k_-)\neq(w_+k_-,k_-)$, then we have $h_{top}(\Omega(k_+,k_-))\neq h_{top}(\Omega(1s,k_-))$.
\end{remark}

We call $\Omega(k_+,k_-)$ is linearizable if there exists linearizable kneading invariants $(1s,0t)$ such that $\Omega(k_+,k_-)=\Omega(1s,0t)$.

\begin{proposition}\label{a=cprop2}
Let $f\in ELM$ with $H=(a,b)$, where $a<c$ is fixed and $b\in E_{f}(a)$. Then we can characterize the maximal plateau $P(b)$ such that for all $\epsilon\in P(b)$, $h_{top}(f|S_{f}(a,\epsilon))=h_{top}(f|S_{f}(a,b))$.
\begin{proof}\rm
Let $(k_{+},k_{-})$ be the kneading invariants of $f$ and $b\in E_{f}(a)$, where $a$ is fixed and $a\neq c$.  By Proposition \ref{bifurcationsetEa}, $b\in E_{f}(a)$ indicates that ${\rm \bf {b}}$ is self-admissible but $({\rm \bf {b}},{\rm \bf {a}})$ may not be weak-admissible. The difference between the case $a=c$ is that ${\rm \bf {b}}$ could also be changed here.
\par  {\bf 1. If ${\rm \bf {b}}$ is periodic.}  By Lemma \ref{exist weak-admissible}, we can always find weak-admissible $(1s,0t)$ such that $\Omega({\rm \bf {b}},{\rm \bf {a}})=\Omega(1s,0t)$. We divide the proof into two main cases.
\par \par Case 1, $\Omega({\rm \bf {b}},{\rm \bf {a}})$ can be non-periodically renormalized during the renormalization process. Suppose $j$-th $(j\geq 1)$ renormalization be the nearest non-periodic renormalization with  words ($w_+,w_-)$, then all the $i$-th ($0\leq i\leq j-1$) renormalizations are periodic with words $(w_{+i},w_{-i})$, and denote $l_i =|w_{+i}|=|w_{-i}|$. By Remark \ref{atleastoncerenor} above,
$$h_{top}(\Omega({\rm \bf {b}},{\rm \bf {a}}))=\frac{\log\beta}{l_1 l_2 \cdots l_{j-1}},$$
where $\log \beta$ denotes the entropy of $(\sigma, \Omega(w_{+}^\infty,w_{-}^\infty))$.
Denote the lexicographically smallest $R^{j-1}(1s)$ be $\eta$ and  the largest $R^{j-1}(1s)$ be $\xi$, as a result, the question here transfers to how to find $\eta$ and $\xi$ via renormalization words $(w_+,w_-)$. By Lemma \ref{w+} above, we obtain that $\xi=w^\infty_+$ and $\eta=w_+ w^\infty_-$.
 Let $P(b)=[b_L,b_R]$, using $\ast$-product,
   $$\tau_{f}(b_{R}+)=(w_{+1},w_{-1})\ast \cdots \ast(w_{+(j-1)},w_{-(j-1)})\ast w_+^\infty, \ \textup{and}$$
   $$\tau_{f}(b_{L}+)=(w_{+1},w_{-1})\ast \cdots \ast(w_{+(j-1)},w_{-(j-1)})\ast w_+w_-^\infty.$$
 We can see from the proof above that the non-periodic renormalization words can directly decide $\eta$ and $\xi$, hence we simply take $j=1$ if $f$ can be non-periodically renormalized, and obtain that $\tau_{f}(b_{R}+)=w_+^\infty$, $\tau_{f}(b_{L}+)=w_+w_-^\infty$. In addition, if we compare the plateau $P(b)$ with the plateau $I(b)$, it is clear that $I(b)$ is only subset of plateau $P(b)$.
\par Case 2, $\Omega({\rm \bf {b}},{\rm \bf {a}})$ can not be non-periodically renormalized, and here we only consider $\Omega({\rm \bf {b}},{\rm \bf {a}})$ being prime. Let $\sigma({\rm \bf {b}})=(v_{1}v_{2}\cdots v_p)^\infty$. Denote $1s=(1v_{1}v_{2}\cdots v_{q-1})^{\infty}$ ($q\leq p$) and $(1s)|_q=(1v_1v_2 \cdots v_{q-1})$. Let $P(b)=[b_L,b_R]$, we have that
$$\tau_{f}(b_{L}+)=(1s)|_q(0t) \ \ {\rm and} \ \ \tau_{f}(b_{R}+)=1s.$$
 Clearly, $I(b)\subsetneqq P(b)$ at this case.
 \par  {\bf 2. If ${\rm \bf {b}}$ is not periodic.}
 We know that $({\rm \bf {b}},{\rm \bf {a}})$ may not be weak-admissible, by Lemma \ref{exist weak-admissible}, there exists weak-admissible $({\rm \bf {b}},0t)$ such that $\Omega({\rm \bf {b}},{\rm \bf {a}})=\Omega({\rm \bf {b}},0t)$. Notice that ${\rm \bf {b}}$ is unchanged here since $I(b)=\{b\}$ by Remark \ref{aperiodic11}. It is clear that $0t$ is periodic if this change happens, and we denote by $0t={\rm \bf {a}}$ if $({\rm \bf {b}},{\rm \bf {a}})$ is weak-admissible. Firstly, if $\Omega({\rm \bf {b}},{\rm \bf {a}})$ can be non-periodically renormalized via $(w_+,w_-)$, including the weak-admissible case $(w_+0t,0t)$, similar to the proof of Case 1 above, $P(b)=[b_L,b_R]$, where $\tau_{f}(b_{R}+)=w_+^\infty$ and $\tau_{f}(b_{L}+)=w_+w_-^\infty$ (and $\tau_{f}(b_{L}+)=w_+0t$ if ${\rm \bf {b}}=w_+0t$). As a result, $I(b)=\{b\}\subsetneqq P(b)$.
 \par Next we consider the case that $\Omega({\rm \bf {b}},{\rm \bf {a}})$ can neither be non-periodically renormalized nor be the case $(w_+0t,0t)$. We prove $P(b)=\{b\}$ by contradiction. Suppose $P(b)\neq \{b\}$, and for any $b^\prime \in P(b)$ with $b\neq b^\prime$, we denote $\tau_f(b^\prime+)={\rm \bf {b}}^\prime$. Naturally, there exists weak-admissible $(1s^\prime,0t^\prime)$ such that $\Omega({\rm \bf {b}}^\prime,{\rm \bf {a}})=\Omega(1s^\prime,0t^\prime)$. Since both $b$ and $b^\prime$ belong to $P(b)$, applying the formula in equation (\ref{gongshi1}), we know that both $\Omega({\rm \bf {b}}^\prime,{\rm \bf {a}})$ and $\Omega({\rm \bf {b}},{\rm \bf {a}})$ correspond to the same $(\beta,\alpha)\in \Delta$.

\par Subcase 1, $\Omega({\rm \bf {b}},{\rm \bf {a}})$ is linearizable. Then there exists $(\beta,\alpha)\in \Delta$ such that $\Omega_{\beta,\alpha}=\Omega({\rm \bf {b}},0t)$.
If $(1s^\prime,0t^\prime)$ is also linearizable, by the proof of Lemma \ref{lemmaxin222}, they correspond to the same $(\beta,\alpha)$, but $1s^\prime\neq{\rm \bf {b}}$, which is impossible. If $(1s^\prime,0t^\prime)$ is not linearizable, $(1s^\prime,0t^\prime)$ could be weak-admissible or non-periodically renormalized, by Lemma \ref{lemmaxin222} and Remark \ref{remarkxin}, $b^\prime \notin P(b)$.
\par Subcase 2, $({\rm \bf {b}},{\rm \bf {a}})$ is weak-admissible and hence $({\rm \bf {b}},{\rm \bf {a}})=({\rm \bf {b}},w_-{\rm \bf {b}})$. We do not consider $({\rm \bf {b}},0t)$ here since $0t$ should be periodic.
This indicates the kneading invariants of $T_{\beta,\alpha}$ is $({\rm \bf {b}},w_-^\infty)$. If $(1s^\prime,0t^\prime)$ is linearizable, then $T_{\beta,\alpha}$ has two different kneading invariants, which is a contradiction. If $(1s^\prime,0t^\prime)$ is not linearizable, then $(1s^\prime,0t^\prime)$ have three cases: 1, $(1s^\prime,0t^\prime)$ can be non-periodically renormalized via $(w_+,w_-)$; 2, $(1s^\prime,0t^\prime)=(1s^\prime,w_-1s^\prime)$; 3, $(1s^\prime,0t^\prime)=(w_+0t^\prime,0t^\prime)$. As a result, we have that the kneading invariants of $T_{\beta,\alpha}$ can be $(w_+^\infty,w_-^\infty)$, $(1s^\prime,w_-^\infty)$ and $(w_+^\infty,0t^\prime)$, respectively. All of them are different from $({\rm \bf {b}},w_-^\infty)$, which leads to the contradiction. $\hfill\square$

\end{proof}
\end{proposition}

\begin{proposition}\label{a=cprop}
Let $f\in ELM$ with $H=(c,b)$, where $b\in E_{f}(c)$. Then we can characterize the endpoints of the maximal plateau $P(b)$ such that for all $\epsilon\in P(b)$, $h_{top}(f|S_{f}(a,\epsilon))=h_{top}(f|S_{f}(a,b))$.
\begin{proof}\rm
Let $(k_{+},k_{-})$ be the kneading invariants of $f$ and $b\in E_{f}(c)$. We know that ${\rm \bf {b}}$ is self-admissible, but $({\rm \bf {b}},k_-)$ may not be weak-admissible since there may exist integer $r$ such that $\sigma^r(k_-)\prec \sigma({\rm \bf {b}})$. By Lemma \ref{exist weak-admissible}, we can always find weak-admissible $({\rm \bf {b}},0t)$ such that $\Omega({\rm \bf {b}},k_-)=\Omega({\rm \bf {b}},0t)$, where $t=(u_1u_2\cdots (u_r-1))^\infty$, and $0t=k_-$ if $r=+\infty$.
\par  {\bf If ${\rm \bf {b}}$ is periodic.}  Let $\sigma({\rm \bf {b}})=(v_{1}v_{2}\cdots v_p)^\infty$.
\par Case 1, $\Omega({\rm \bf {b}},k_-)$ can be non-periodically renormalized via words ($w_+,w_-)$, including the case $(w_+^\infty,w_-w_+^\infty)$. Similar to Case 1 in the proof of Proposition \ref{a=cprop2}, we have that $P(b)=[b_L,b_R]$, where
$ \tau_{f}(b_{L}+)=w_+w_-^\infty \ \ {\rm and} \ \ \tau_{f}(b_{R}+)=w_+^\infty.$ It is clear that $I(b)\subsetneqq P(b)$.

\par Case 2, $({\rm \bf {b}},k_-)$ is linearizable ($r=+\infty$). Here we only consider $f$ being prime. Denote ${\rm \bf {b}}|_p=(1v_1v_2 \cdots v_{p-1})$. Let $P(b)=[b_L,b_R]$, by Lemma \ref{normalrenor},
$\tau_{f}(b_{L}+)={\rm \bf {b}}|_p1\sigma(k_+)={\rm \bf {b}}|_pk_+ \ \ {\rm and} \ \ \tau_{f}(b_{R}+)={\rm \bf {b}}.$ Notice that $\tau_{f}(b_{L}-)={\rm \bf {b}}|_pk_-$, hence $P(b)=I(b)=[b_L,b_R]$.
\par Case 3, $({\rm \bf {b}},0t)$ is linearizable ($r<+\infty$). By Case 2, we have that $P(b)=[b_L,b_R]$, where $\tau_{f}(b_{L}+)={\rm \bf {b}}|_p(0t)\prec1v_{1}v_{2}\cdots v_{p-1}k_- \ \ {\rm and} \ \ \tau_{f}(b_{R}+)={\rm \bf {b}}.$ By Theorem \ref{plateauxin}, $I(b)=[b_L^\prime,b_R^\prime]$, where $\tau_{f}(b_{L}^\prime+)={\rm \bf {b}}|_pk_+$ and $\tau_{f}(b_{R}^\prime+)={\rm \bf {b}}$. Hence $b_L < b_L^\prime$
 and $I(b)\subsetneqq P(b)$.
 \par  {\bf If ${\rm \bf {b}}$ is not periodic.} Firstly, if $\Omega({\rm \bf {b}},k_-)$ can be non-periodically renormalized via $(w_+,w_-)$, similar to Case 1 above, $P(b)=[b_L,b_R]$, where $\tau_{f}(b_{R}+)=w_+^\infty$ and $\tau_{f}(b_{L}+)=w_+w_-^\infty$. Hence $I(b)=\{b\}\subsetneqq P(b)$ at this case. Notice that we do not need to consider the weak-admissible cases $(w_+k_-,k_-)$ and $(w_+0t,0t)$  here, since it can be concluded by Case 2 and Case 3 above with ${\rm \bf {b}}$ being periodic. Similarly to the proof of Proposition \ref{a=cprop2}, when considering the cases that $\Omega({\rm \bf {b}},k_-)$ can not be non-periodically renormalized or be the weak-admissible case $({\rm \bf {b}},w_-{\rm \bf {b}})$, we have $P(b)=I(b)=\{b\}$. $\hfill\square$

\end{proof}
\end{proposition}

As we can see in Example \ref{different}, if we consider $f\in ELM$ with $(k_+,k_-)=(10^\infty,(01110)^\infty)$, and the hole $(c,b)$. We can see that $f,g,h$ in Example \ref{different} correspond to three different bifurcation points in $E_f(c)$, but they are on the plateau $P(b)$ of identical entropy $\log1.3247$. Especially, when $\tau_f(b+)=10(011)^\infty$ which is not periodic, $I(b)=\{b\}$ but $P(b)\neq\{b\}$.

\vspace{0.2cm}
\noindent {\bf Proof of Theorem \ref{plateau}}
\par Let $f\in ELM$ with $H=(a,b)$, where $a\in (0,c]$ is fixed and $b\in E_{f}(a)$. With the help of Proposition \ref{a=cprop2} and Proposition \ref{a=cprop}, we are able to characterize the endpoints of the maximal plateau $P(b)$ such that for all $\epsilon\in P(b)$, $h_{top}(f|S_{f}(a,\epsilon))=h_{top}(f|S_{f}(a,b))$. Moreover, at some special cases, $I(b)=P(b)=\{b\}$ is a singleton, see the following two remarks.
 $\hfill\square$

\begin{example}\label{consturcitionbLbR22}(Plateau of $h_{top}(\tilde{S}_{f}(a,b))$)\rm\begin{enumerate}\item Let $f$ be the doubling map and $(a,b)$ be the hole, where ${\rm \bf {b}}=(10011)^\infty$ and ${\rm \bf {a}}=(0111010100)^\infty$. We can see that $\Omega({\rm \bf {b}},{\rm \bf {a}})=\Omega((10011)^\infty,(0111010)^\infty)$ and $({\rm \bf {b}},{\rm \bf {a}})$ can be renormalized via $(w_+,w_-)=(10,011)$. We have that $P(b)=[b_{L},b_{R}]$, where $\tau_{f}(b_{L}+)=10(011)^\infty$ and $\tau_{f}(b_{R}+)=(10)^\infty$.  Moreover, $I(b)=(b_l,b_r]$ where $\tau_{f}(b_{l}+)=10011(0111010)^\infty$ and $\tau_{f}(b_{r}+)=(10011)^\infty$.  Hence $\tau_{f}(b_{L}+)\prec \tau_{f}(b_{l}+)$ and $\tau_{f}(b_{r}+)\prec \tau_{f}(b_{R}+)$, $I(b)\subsetneqq P(b)$.\item Let $f\in ELM$ with $(k_+,k_-)=((10000)^\infty,(011)^\infty)$ and $(a,b)$ be the hole, where ${\rm \bf {b}}=(100)^\infty$ and ${\rm \bf {a}}=(011000)^\infty$. By Lemma \ref{exist weak-admissible}, $\Omega({\rm \bf {b}},{\rm \bf {a}})=\Omega((100)^\infty,(01)^\infty)$ and hence $({\rm \bf {b}},{\rm \bf {a}})$ is prime. By the proof of Subcase 2 in Theorem \ref{plateau}, we have that $P(b)=[b_{L},b_{R}]$, where $\tau_{f}(b_{L}+)=100(01)^\infty$ and $\tau_{f}(b_{R}+)=(100)^\infty$. For any $\gamma\in P(b)$, $h_{top}(\tilde{S}_{f}(a,\gamma))=h_{top}(\Omega((100)^\infty,(01)^\infty))$. Moreover, $I(b)=(b_L,b_R]$ and$b_L\notin I(b)$ for the reason that $\Omega(100(01)^\infty,(01)^\infty)\neq \Omega((100)^\infty,(01)^\infty)$, and $b_L\in E_{f}(a)$. Hence $I(b)\subsetneqq P(b)$.\item Also let $f\in ELM$ with $(k_+,k_-)=((10000)^\infty,(011)^\infty)$ and $(a,b)$ be the hole, while ${\rm \bf {b}}=(100)^\infty$ and ${\rm \bf {a}}=(01100)^\infty$. Similarly, $\Omega({\rm \bf {b}},{\rm \bf {a}})=\Omega((100)^\infty,(01)^\infty)$, $({\rm \bf {b}},{\rm \bf {a}})$ is prime, and $P(b)=[b_{L},b_{R}]$, where $\tau_{f}(b_{L}+)=100(01)^\infty$ and $\tau_{f}(b_{R}+)=(100)^\infty$. However, different from (ii) above, here $I(b)=(b_l,b_r]$, where $\tau_{f}(b_{l}+)=100(0110)^\infty$ and $\tau_{f}(b_{r}+)=(100)^\infty$.Hence $I(b)\subsetneqq P(b)$ for $\tau_{f}(b_{L}+)\prec\tau_{f}(b_{l}+)$.\end{enumerate}\end{example}

\begin{remark}\label{anot=c}

Let $b\in E_f(a)$ and ${\rm \bf {b}}$ is periodic.
 \begin{enumerate}
\item The plateau $P(b)$ is always a closed subinterval. When $a\neq c$, $I(b)\subsetneqq P(b)$ for $I(b)$ is always left open; when $a=c$, $I(b)=P(b)$ only when $({\rm \bf {b}},k_-)$ is linearizable.
\item We can see that $b\in P(b)$ but $b$ may not be the right endpoint of $P(b)$, which is quite different from Remark \ref{aperiodic11}, see Example \ref{different}.
 \end{enumerate}
\end{remark}

\begin{remark}\label{xinremarka}
Let $b\in E_{f}(a)$ and ${\rm \bf {b}}$ is not periodic. Denote by ${\rm \bf {a}}=k_-$ at the case $a=c$. Then $P(b)=I(b)=\{b\}$ if and only if $\Omega({\rm \bf {b}},{\rm \bf {a}})$ is linearizable or $({\rm \bf {b}},{\rm \bf {a}})=({\rm \bf {b}},w_-{\rm \bf {b}})$.
\end{remark}

\section{Two bifurcation sets coincide}\label{Some examples}
It was proved in {\cite[Theorem 1.6]{sun2024}} that $E_{f}(a)$ is of null Lebesgue measure, and hence the entropy function $\lambda_{f}(a):b\mapsto h_{top}(f|S_{f}(a,b))$ is a devil staircase for the reason that $B_{f}(a)\subseteq E_{f}(a)$. A natural question arises: when will the two bifurcation sets coincide? We answer the question in two cases: the case $a=c$ and the case $a\neq c$.

\begin{proposition}\label{lemma28}
Let $a\in (0,c]$ being fixed and $b\in E_f(a)$ with ${\rm \bf {b}}$ not being periodic. Then there exists $b^\prime\in E_f(a)$ with ${\rm \bf {b}}^\prime$ being periodic such that $d({\rm \bf {b}},{\rm \bf {b}}^\prime)<\epsilon$ for any $\epsilon>0$.
\begin{proof} \rm Denote ${\rm \bf {b}}=(v_{1}v_{2}\cdots )$ and ${\rm \bf {a}}=( u_{1}u_{2}\cdots )$. By Remark \ref{aperiodic11}, $I(b)=\{b\}$, hence ${\rm \bf {b}}$ will not be changed and ${\rm \bf {a}}$ maybe changed into a periodic sequence. Since ${\rm \bf {b}}\neq 10^\infty$, choose any $v_n=1$ and there exists a minimal integer $j\leq n-1$ such that $(v_1\cdots v_{n-j})=(v_{j+1}\cdots v_n)$, and the existence of such $j$ is secured by $j=n-1$. The integer $n$ such that $v_n=1$ can be arbitrarily large. Since ${\rm \bf {b}}$ is self-admissible, we have $\sigma^j({\rm \bf {b}})\succ {\rm \bf {b}}$ and there exists a minimal integer $r\geq n$ such that $(v_1\cdots v_{n-j}\cdots v_{r-j})=(v_{j+1}\cdots v_n \cdots v_r)$, $v_{r+1}=1$ and $v_{r-j+1}=0$. Notice that $\sigma^j({\rm \bf {b}})\neq {\rm \bf {b}}$ for ${\rm \bf {b}}$ not being periodic. Let ${\rm \bf {b}}^\prime=(v_1\cdots v_{n-j}\cdots v_{r})^\infty$, next we prove that ${\rm \bf {b}}^\prime$ is self-admissible, that is, $\sigma({\rm \bf {b}}^\prime)\preceq \sigma^{m}({\rm \bf {b}}^\prime)$ for all $m\in\{0,\cdots,r-1\}$. If $v_{m+1}=1$, $\sigma({\rm \bf {b}}^\prime)\preceq \sigma^{m}({\rm \bf {b}}^\prime)$ is obvious. It remains to show that ${\rm \bf {b}}^\prime\preceq\sigma^m({\rm \bf {b}}^\prime)$ for all $m$ such that $v_{m+1}=1$, as this implies that $\sigma({\rm \bf {b}}^\prime)\preceq \sigma^{m}({\rm \bf {b}}^\prime)$ when $v_{m+1}=0$. Let $v_{m+1}=1$, if $1\leq m< j$, then ${\rm \bf {b}}\preceq\sigma^m({\rm \bf {b}})$ and the minimality of $j$ implies that $(v_{m+1},\cdots, v_n)\succ (v_1,\cdots,v_{n-m})$, hence ${\rm \bf {b}}^\prime\prec\sigma^m({\rm \bf {b}}^\prime)$. If $m\geq j$, then we have
\begin{eqnarray*}
(v_{m+1}\cdots v_{r} v_{1})&=(v_{m+1}\cdots v_{r} v_{r+1})\!\\
 &\succ(v_{m+1}\cdots v_{r} 0)\!\\
 &=(v_{m-j+1}\cdots v_{r-j} v_{r-j+1})\!\\
  &\succeq(v_{1}\cdots v_{r-m} w_{r-m+1})\!
\end{eqnarray*}
when $v_{m+1}=1$, which yields ${\rm \bf {b}}^\prime\prec\sigma^m({\rm \bf {b}}^\prime)$. As a result, ${\rm \bf {b}}^\prime\in\Omega({\rm \bf {b}}^\prime,{\rm \bf {a}})$, ${\rm \bf {b}}^\prime$ is self-admissible and $b^\prime\in E_f(a)$.

\par As we can see, ${\rm \bf {b}}^\prime\prec{\rm \bf {b}}$ makes sure that ${\rm \bf {b}}^\prime$ approaches ${\rm \bf {b}}$ from the left side.
Now we consider the distance $d({\rm \bf {b}}^\prime,{\rm \bf {b}})$, where $d(,)$ is the usual metric on $\{0,1\}^{\mathbb{N}}$. Since  the integer $n$ can be arbitrarily large, $r\geq n$ is also arbitrarily large. Then the distance $d({\rm \bf {b}}^\prime,{\rm \bf {b}})=2^{-r}$ can be arbitrarily small and hence $d({\rm \bf {b}}^\prime,{\rm \bf {b}})<\epsilon$ for any $\epsilon>0$. $\hfill\square$
\end{proof}
\end{proposition}

\vspace{0.2cm}
\noindent {\bf Proof of Theorem \ref{thmxin}}
\par Let $f\in ELM$ with $H=(a,b)$, where $a\in (0,c]$ is fixed and $b\in E_{f}(a)$. Denote by $ D_f(a):=\{b\in E_f(a):{\rm \bf {b}} \ \textup {is periodic}\}$. By the proof of Proposition \ref{lemma28}, for any $b\in E_f(a)$ with ${\rm \bf {b}}$ not being periodic, we can always find $b^\prime\in E_f(a)$ with ${\rm \bf {b}}^\prime$ being periodic, and the distance between ${\rm \bf {b}}$ and ${\rm \bf {b}}^\prime$ is arbitrarily small. Since $f\in ELM$ is topologically expansive, we have that the Euclidean distance between two points $b$ and $b^\prime$ can also be arbitrarily small. Hence the set $D_f(a)$ is dense in $E_f(a)$. The proof is completed. $\hfill\square$
\vspace{0.2cm}
\begin{lemma}\label{abudengyuc}
Let $a\neq c$. Then $B_{f}(a)= E_{f}(a)$ if and only if for all $b\in E_{f}(a)$, ${\rm \bf {b}}$ satisfies
 \begin{enumerate}
\item ${\rm \bf {b}}$ is not periodic,
\item $\Omega({\rm \bf {b}},{\rm \bf {a}})$ is linearizable or $({\rm \bf {b}},{\rm \bf {a}})=({\rm \bf {b}},w_-{\rm \bf {b}})$.
\end{enumerate}

\begin{proof}\rm
\par  Let $b\in E_f (a)$. If ${\rm \bf {b}}$ is periodic, by by Remark \ref{anot=c} (i), $P(b)$ is always closed while $I(b)$ is always left open, hence $I(b)\subsetneqq P(b)$ and $E_f(a)\neq B_{f}(a)$. If ${\rm \bf {b}}$ is not periodic and $\Omega({\rm \bf {b}},{\rm \bf {a}})$ is not linearizable, by the proof of Proposition \ref{a=cprop2}, $I(b)=\{b\}$ while $P(b)$ is a proper subinterval, hence we still have $I(b)\subsetneqq P(b)$ and $E_f(a)\neq B_{f}(a)$. For the case that ${\rm \bf {b}}$ is not periodic and $\Omega({\rm \bf {b}},{\rm \bf {a}})$ can be linearizable or $({\rm \bf {b}},{\rm \bf {a}})=({\rm \bf {b}},w_-{\rm \bf {b}})$, also by Proposition \ref{a=cprop2}, $P(b)=I(b)=\{b\}$.
$\hfill\square$
\end{proof}
\end{lemma}

\begin{remark}\label{remark25}
Let $a\neq c$ and $b\in E_f(a)$ with ${\rm \bf {b}}$ not being periodic. By Proposition \ref{lemma28} and Lemma \ref{abudengyuc}, there always exists $b^\prime\in E_f(a)$ with ${\rm \bf {b}}^\prime$ being periodic. By Proposition \ref{a=cprop2}, we have $I(b^\prime)\subsetneqq P(b^\prime)$ and immediately $E_f(a)\neq B_f(a)$. As a result, $E_f(a)\neq B_f(a)$ when $a\neq c$.

\end{remark}

\begin{lemma}\label{lemma26}
Let $b\in E_f(a)$ with $a\in (0,c]$ being fixed, where ${\rm \bf {b}}$ is not periodic, $\Omega({\rm \bf {b}},{\rm \bf {a}})$ is linearizable and $({\rm \bf {b}},{\rm \bf {a}})$ is not admissible. Then there exists $b^\prime\in E_f(a)$ with periodic ${\rm \bf {b}}^\prime$ such that $({\rm \bf {b}}^\prime,{\rm \bf {a}})$ is not admissible.

\begin{proof} \rm
Denote ${\rm \bf {b}}=(v_{1}v_{2}\cdots )$. The assumption $({\rm \bf {b}},{\rm \bf {a}})$ not being admissible indicates that there exists a minimal integer $i$ such that $\sigma^{i}({\rm \bf {a}})\prec {\rm \bf {b}}$. Hence we can always find an integer $n>i$ such that two finite words $(v_{1}v_{2}\cdots v_n)=(u_{i+1}\cdots u_{i+n})$ and $v_{n+1}=1$, $u_{i+n+1}=0$. Without loss of generality, we only consider the case $v_n=u_{i+n}=1$, and the case $v_n=u_{i+n}=0$ can be proved similarly. If $(v_{1}\cdots v_n)^\infty$ is self-admissible, then we let ${\rm \bf {b}}^\prime=(v_{1}\cdots v_n)^\infty$, it is clear that $b^\prime\in E_f(a)$ and $({\rm \bf {b}}^\prime,{\rm \bf {a}})$ is not admissible. If $(v_{1}\cdots v_n)^\infty$ is not self-admissible, then there exists a minimal integer $j\leq n-1$ such that $(v_1\cdots v_{n-j})=(v_{j+1}\cdots v_n)$. Since ${\rm \bf {b}}$ is self-admissible, we have $\sigma^j({\rm \bf {b}})\succ {\rm \bf {b}}$ and there exists a minimal integer $r\geq n$ such that $(v_1\cdots v_{n-j}\cdots v_{r-j})=(v_{j+1}\cdots v_n \cdots v_r)$, $v_{r+1}=1$ and $v_{r-j+1}=0$. Let ${\rm \bf {b}}^\prime=(v_{1}\cdots v_{r-j})^\infty$, we can see that $\sigma^{i}({\rm \bf {a}})\prec{\rm \bf {b}}\prec{\rm \bf {b}}^\prime\prec\sigma^{j}({\rm \bf {b}})$. Hence ${\rm \bf {b}}^\prime$ is self-admissible and $({\rm \bf {b}}^\prime,{\rm \bf {a}})$ is not admissible.
$\hfill\square$
\end{proof}
\end{lemma}

\begin{lemma}\label{lemma27}
Let $b\in E_f(a)$ with ${\rm \bf {b}}$ not being periodic and ${\rm \bf {a}}=w_-{\rm \bf {b}}$, where $({\rm \bf {b}},{\rm \bf {a}})$ is weak-admissible. Then there exists $b^\prime\in E_f(a)$ with periodic ${\rm \bf {b}}^\prime$ such that $({\rm \bf {b}}^\prime,{\rm \bf {a}})$ is not admissible.

\begin{proof} \rm
Denote the length of word $w_-$ as $p$ and ${\rm \bf {b}}=(v_{1}v_{2}\cdots )$.
Similar to the proof of Lemma \ref{lemma26}, choose any $v_n=1$ and there exists an integer $j$ such that $(v_1\cdots v_{n-j})=(v_{j+1}\cdots v_n)$. Since ${\rm \bf {b}}$ is self-admissible, we have $\sigma^j({\rm \bf {b}})\succ {\rm \bf {b}}$ and there exists a minimal integer $r\geq n$ such that $(v_1\cdots v_{n-j}\cdots v_{r-j})=(v_{j+1}\cdots v_n \cdots v_r)$, $v_{r+1}=1$ and $v_{r-j+1}=0$. Let ${\rm \bf {b}}^\prime=(v_1\cdots v_{r-j})^\infty$, it can be verified that ${\rm \bf {b}}^\prime$ is self-admissible, which indicates $b^\prime\in E_f(a)$. Moreover, $\sigma^{p}({\rm \bf {a}})={\rm \bf {b}}\prec{\rm \bf {b}}^\prime\prec\sigma^{j}({\rm \bf {b}})$, hence $({\rm \bf {b}}^\prime,k_-)$ is not admissible. $\hfill\square$

\end{proof}
\end{lemma}

\begin{remark}\label{remark26}
By Lemma \ref{lemma27}, if $k_-$ is not periodic and can be written into the form $k_-=w_-{\rm \bf {b}}$, where both $w_-$ and ${\rm \bf {b}}$ are self-admissible. Then $E_f(c)\neq B_f(c)$.

\end{remark}

\vspace{0.2cm}

\noindent {\bf Proof of Theorem \ref{thmplat}}
\par Let $f\in ELM$ with a hole $H=(a,b)$, and $(k_{+},k_{-})$ be its kneading invariants. For the case $a\neq c$, by Remark \ref{remark25}, we can obtain that $E_f(a)\neq B_f(a)$. Next we focus on the case $a=c$.
\par Let $b\in E_f(c)$. When $k_-$ is periodic, by Remark \ref{anot=c}, we have that $I(b)=P(b)$ if and only if $({\rm \bf {b}},k_-)$ is linearizable. When $k_-$ is not periodic, by Remark \ref{xinremarka}, $I(b)=P(b)=\{b\}$ if and only if $\Omega({\rm \bf {b}},{\rm \bf {a}})$ is linearizable or $({\rm \bf {b}},k_-)=({\rm \bf {b}},w_-{\rm \bf {b}})$. Applying the proof of Lemma \ref{lemma26} and \ref{lemma27}, if $\Omega({\rm \bf {b}},{\rm \bf {a}})$ is linearizable or $({\rm \bf {b}},k_-)=({\rm \bf {b}},w_-{\rm \bf {b}})$, then there exists $b^\prime\in E_f(c)$ with periodic ${\rm \bf {b}}^\prime$ such that $({\rm \bf {b}}^\prime,{\rm \bf {a}})$ is not linearizable. As we know, $E_f(c)=B_f(c)$ is equivalent to that $I(b)=P(b)$ for any $b\in E_f(c)$. Hence at the case $a=c$, $E_f(c)=B_f(c)$ if and only if for all $b\in E_{f}(c)$, $({\rm \bf {b}},k_-)$ is linearizable.  $\hfill\square$

\section*{Data Availability}
 All data generated or analysed during this study are included in this published article.

\section*{Conflict of interest}
 The authors declare that they have no known
competing financial interests or personal relationships that could have appeared to influence the work reported
in this paper.

\section*{Acknowledgment}
Y. Sun was supported by CSC 202306150094, B. Li was supported by NSFC 12271176 and Guangdong Natural Science Foundation 2023A1515010691. We appreciate the anonymous referees whose constructive comments helped to improve the article significantly.

\end{document}